\documentclass[preprint,11pt]{elsarticle}
\usepackage{mathrsfs}
\usepackage{amsmath}
\usepackage{amsthm}
\usepackage{amssymb}
\usepackage{graphicx}
\usepackage{color}
\usepackage{amsfonts}
\usepackage{lineno}

\newtheorem{theorem}{Theorem}[section]
\newtheorem{lemma}[theorem]{\textbf{Lemma}}
\newtheorem{corollary}[theorem]{\textbf{Corollary}}
\newtheorem{proposition}[theorem]{\textbf{Proposition}}

\numberwithin{equation}{section}
\journal{Frontiers of Mathematics in China}

\begin{document}

\begin{frontmatter}

\title{The contact process on the regular tree with random vertex weights}

\renewcommand{\thefootnote}{\fnsymbol{footnote}}
\author{Pan Yu\footnote{\textbf{E-mail}: perryfly@pku.edu.cn \textbf{Address}: School of Mathematical Sciences, Peking University, Beijing 100871, China.}, Chen Dayue}
\address{LMAM, Peking University}
\author{Xue Xiaofeng}
\address{School of Science, Beijing Jiaotong University}
\begin{abstract}
This paper is concerned with contact process with random vertex weights on regular trees, and study the asymptotic behavior of the critical infection rate as the degree of the trees increasing to infinity. In this model, the infection propagates through the edge connecting vertices $x$ and $y$ at rate $\lambda\rho(x)\rho(y)$ for some $\lambda>0$, where $\{\rho(x),{x\in T^d}\}$ are $i.i.d.$ vertex weights. We show that when $d$ is large enough there is a phase transition at $\lambda_c(d)\in(0,\infty)$ such that for $\lambda<\lambda_c(d)$ the contact process dies out, and for $\lambda>\lambda_c(d)$ the contact process survives with a positive probability. Moreover, we also show that there is another phase transition at $\lambda_e(d)$ such that for $\lambda<\lambda_e(d)$ the contact process dies out at an exponential rate. Finally, we show that these two critical values have the same asymptotic behavior as $d$ increases.
\end{abstract}
\end{frontmatter}

\section{Introduction}
We first give the description of the process (denoted by $\{\eta_t\}_{t\geq0}$) studied in this paper. Let $T^d$ be a rooted regular tree where the root $O$ has degree $d$ and each of other vertices has degree $d+1$. On each vertex $x\in T^d,$ we put
a nonnegative random variable denoted by $\rho(x).$ We call this nonnegative random variable $\rho(x)$ the random weight on $x.$ We omit the superscript or subscript $d$ when there is no misunderstanding.

The state space of the process is $X=\{0,1\}^{T}$ with the product topology. For given $\{\rho(x),{x\in T^d}\}$, the process is a particle system described by a collection of transition measures:
\[c(x,\eta)=\left\{ \begin{array}{ll}
1,& \ {\rm if}\ \eta(x)=1,\\
\lambda\sum_{y:y\sim x}\rho(x)\rho(y)\eta(y),& \ {\rm{if}}\ \eta(x)=0,
\end{array} \right.\]
where $\lambda$ is a nonnegative parameter.

One interpretation of this process is as a model for the spread of an infection. An individual at $x\in T$ is infected if $\eta(x)=1$ and healthy if $\eta(x)=0$. A healthy individual $x$ will be infected by some infected neighbor $y$ at rate $\lambda\rho(x)\rho(y).$ Infected individuals recover at a constant rate, which is normalized to be 1. We will make the following assumption on the vertex weights.

\textbf{Assumption :} The sequence of vertex weights $\{\rho(x),{x\in T}\}$ is an $i.i.d.$  sequence of nonnegative random variables. Moreover, each $\rho(x)$ is a bounded (by some positive number $M$) random variable.

The boundedness condition ensures the process is a spin system.(See \textbf{Appendix} \textit{5.1}  for more details)

The classic contact process on integer lattice $\mathbb{Z}^d$ was introduced by Harris \cite{1974}. In \cite{1974}, Harris shows that the process exhibits a phase transition as the parameter $\lambda$ changes. In \cite{1}, Harris introduces the graphical representation, that is a useful method in studying the contact process.  In \cite{bg}, Bezuidenhout and Grimmett give the whole results about the convergence theorem of the process. They show that the critical contact process on $\mathbb{Z}^d$ dies out and the complete convergence theorem  holds for contact process on $\mathbb{Z}^d$.

For the model in this paper, let the vertex weights $\rho\equiv1$, then it is the contact process on trees. In \cite{1992}, Pemantle introduces the contact process on trees, which exhibits a different property from the classic contact process. For readers who want to learn more about the contact process, the books \cite{ips} and \cite{4} is a good reference.

Recently, the contact process in random environment has been a direction of intensive study. In \cite{bds}, Bramson, Durrett and Schonmann study the contact process in a random environment, where the recovery rates are random variables. They show that in one dimension, this process has an intermediate phase in which it survives but not grow linearly. In \cite{liggett1990}, Liggett studies the contact process with random infection rates and recovery rates on $\mathbb{Z}$. Liggett gives sufficient conditions for extinction of the process, for survival of the process and for the process to have at most four extremal invariant measures.  More related results can be found in \cite{Andjel}, \cite{klein}, \cite{liggett1992} and \cite{newman}. In \cite{yc}, Chen and Yao study the contact process on the open clusters of bond percolation on lattices. They show the complete convergence theorem holds for their model. In \cite{r}, Remenik study the contact process in a dynamic random environment, where the environment changes throughout time. Remenik shows the critical process dies out and that the complete convergence theorem holds in the supercritical case.

If we let the vertex weights $\rho$ be binomial, then the model in this paper is a contact process on Galton-Watson trees. In \cite{2001}, Pemantle and Stacey study contact processes on Galton-Watson trees. They show that the contact process has only one phase transition on trees satisfying some conditions. In \cite{2011}, Peterson introduces the contact processes with random vertex weights on complete graphs. Peterson shows that the critical value of this model is the inverse of the second moment of the vertex weight. In \cite{2013}, \cite{2015} and \cite{ar}, it is also shown that the critical value of the model is inversely proportional to the second moment of the vertex weight.

\section{Main results}
Before stating our main results, we first give some notations and definitions. The random vertex weights $\{\rho(x),{x\in T}\}$ are defined on some probability space denoted by $(\Omega,\mathcal{F},\mu)$. Let $E$ be the expectation operator corresponding to the measure $\mu$. $\forall\omega \in \Omega$, take $P_{\lambda}^{\omega}$ as the law of the contact process on $T$ with infection rate $\lambda$ and vertex weights $\{\rho(x,\omega),{x\in T}\}$. The expectation operator corresponding to $P_{\lambda}^{\omega}$ is denoted by $E_{\lambda}^{\omega}$. Define the annealed measure
\[\mathbb{P}_{\lambda}(\cdot)=EP_{\lambda}^{\omega}(\cdot)=\int P_{\lambda}^{\omega}(\cdot)\mu(d\omega).\]
The expectation operator corresponding to $\mathbb{P}_{\lambda}$ is denoted by $\mathbb{E}_{\lambda}$. Take $\mathbb{P}_{\lambda}^d$  (or $\mathbb{E}_{\lambda}^d$) as $\mathbb{P}_{\lambda}$ (or $\mathbb{E}_{\lambda}$) when we need to emphasize the degree $d$ of the tree $T$.

For $t\geq0,$ let
$C_t=\{x\in T:\eta_t(x)=1\}$ be the set of the infected individuals at time $t.$ Assume $C_0=\{O\}$ in this section. Since $\{\eta_t\}_{t\geq0}$ is attractive in the sense of section 3.2 of \cite{ips}, we get the following monotonicity by the basic coupling (see section 3.2 of \cite{ips}):
\begin{align}\label{2.1}
\mathbb{P}_{\lambda_1}(C_t\neq\varnothing,\forall t\geq0)\leq \mathbb{P}_{\lambda_2}(C_t\neq\varnothing,\forall t\geq0)
\end{align}
whenever $\lambda_1<\lambda_2.$

Hence, by \eqref{2.1}, define the critical value
\[\lambda_c(d)=\sup\{\lambda:\mathbb{P}_{\lambda}^d(C_t\neq\varnothing,\forall t\geq0)=0\}.\]

Thus, if $\lambda<\lambda_c,$ we have, $\mathbb{P}_\lambda^d(C_t\neq\varnothing)\rightarrow0,$ as $t\rightarrow\infty.$ So we care about the convergence rate. Therefore, we define another critical value:
\[\lambda_e(d)=\sup\{\lambda:\limsup_{t\rightarrow\infty}\frac{1}{t}\log \mathbb{P}_\lambda^d(C_t\neq\varnothing)<0\}.\]

We call $\lambda_e$ the exponential critical value. It is easy to show that $\lambda_e\leq\lambda_c$.
Note that we do not even know whether $\lambda_c\in(0,\infty)$ or not.

\begin{theorem}\label{thm1}
Assume that $\mu(\rho>0)>0$, then
\begin{align}\label{2.2}
\lambda_c(d)\geq\lambda_e(d)\geq(dE{\rho}^2+\frac{M^4}{E{\rho}^2})^{-1}>0.
\end{align}

Furthermore, if $\frac{(1+\lambda M^2)^2}{\lambda E\rho^2}<d$ for some $\lambda\in(0,\infty)$, then
\[\lambda_e\leq\lambda_c\leq\lambda<\infty.\]
\end{theorem}

Note that the existence of $\lambda\in(0,\infty)$ such that $\frac{(1+\lambda M^2)^2}{\lambda E\rho^2}<d$ can be easily checked since it is a condition of
a quadratic inequality. Anyway, this condition can be satisfied when $d$ is large enough.

\begin{corollary}\label{thm main}
Assume that $\mu(\rho>0)>0$,\ then
\[\lim_{d\rightarrow\infty}d\lambda_c(d)=\lim_{d\rightarrow\infty}d\lambda_e(d)=\frac{1}{E{\rho}^2}.\]
\end{corollary}

\textbf{Theorem 2.1} shows that there is a critical value of $\lambda$ for which $\{\eta_t\}_{t\geq0}$ exhibits a phase transition under some conditions. Since it is difficult to calculate the exact value of the critical value, we would like to study the asymptotic behavior of the critical value. \textbf{Corollary 2.2} shows that $\lambda_c(d)\approx1/(dE\rho^2)$ when $d$ is large enough, which is inversely proportional to the degree
of the root and the second moment of the vertex weight. In the case of  classical contact process on $\mathbb{Z}^d$, the critical value $\lambda_c(d)\approx1/(2d)$ which satisfies the similar law. In \cite{2011}, it is also shown that the critical value of the model is inversely proportional to the second moment of the vertex weight. Related examples can also be found in \cite{2013}, \cite{2015} and \cite{ar}.

When $\rho\equiv1,$ \textbf{Theorem 2.1} and \textbf{Corollary 2.2} show that $\lambda_c\geq1/(d+1)$ and
\[\lim_{d\rightarrow\infty}d\lambda_c(d)=\lim_{d\rightarrow\infty}d\lambda_e(d)=1,\]
that is the similar results with the contact process on regular tree shown by Pemantle \cite{1992}.

\textbf{Proof of Corollary 2.2.} By \eqref{2.2},
\begin{align}\label{lowlim}
\liminf_{d\rightarrow\infty}d\lambda_c(d)\geq\liminf_{d\rightarrow\infty}d\lambda_e(d)\geq\frac{1}{E{\rho}^2}.
\end{align}

Now take $\gamma>1$ and $\lambda=\frac{\gamma}{dE\rho^2},$ then
\[\frac{(1+\lambda M^2)^2}{d\lambda E\rho^2}=\frac{1}{\gamma}(1+\frac{\gamma}{dE\rho^2}M^2)^2\rightarrow \frac{1}{\gamma}<1,\ d\rightarrow \infty.\]

So $\lambda_c(d)\leq \lambda=\frac{\gamma}{dE\rho^2}$ when $d$ is large enough by \textbf{Theorem 2.1}.

Hence,
\begin{align}\label{2.4}
\limsup_{d\rightarrow\infty}d\lambda_c(d)\leq\frac{\gamma}{E\rho^2}.
\end{align}

Let $\gamma\downarrow1,$  the corollary follows from \eqref{lowlim} and \eqref{2.4}.$~~~~~~\Box$

$Some~words~about~the~proof~of~\textbf{Theorem 2.1}:$ The proof is split in section 3 and section 4. First, we show the lower bound of $\lambda_e$ in section 3. There are two key points in the proof. The first one is inspired by \cite{bcp} which revealed the connection between contact process and binary contact path process. This connection helps us convert estimating survival probability into solving  a linear ODE. In this part, the main difficulty is how to handle with the expectation of the product of the vertex weights in a path with given length and may backtrack (see section 3 for more details). There, we introduce an auxiliary random walk which reduces the difficulty of estimates. This is motivated by Kesten in \cite{kesten}. In section 4, we introduce the SIR epidemic model to dominate the process $\{\eta_t\}_{t\geq0}$, and then get an upper bound of $\lambda_c$ by the second-moment method and some sophisticated estimates.

\section{Lower bound for $\lambda_e$}
\label{S:3}
In this section, we will get a lower bound of $\lambda_e$ to show \eqref{2.2}. Throughout this section, we assume the contact process $\{{\eta}_{t}\}_{t\geq0}$ has initial state $\eta_0\equiv1.$

We first introduce a process $\{\xi_{t}\}_{t\geq0}$ which has a close relationship to the process $\{\eta_t\}_{t\geq0}$. The idea of the introduction of the process $\{\xi_{t}\}_{t\geq0}$ comes from Griffeath in \cite{bcp}. In this paper, $\{\xi_{t}\}_{t\geq0}$ is a linear system with values in $\{0,1,2,\cdots\}^{T}$ and the evolution is described as follows when the random environment $\{\rho(x),{x\in T}\}$ is given.

For each $x\in T,$ $\xi(x)\rightarrow 0$ at rate 1. For each nearest neighbor pair $x$ and $y$, $\xi(y)\rightarrow\xi(x)+\xi(y)$ at rate $\lambda\rho(x)\rho(y).$ We also assume $\xi_0\equiv1$.

Formally, the generator of this linear system is
\[
\mathcal{A}f(\xi)=\sum_{x:x\in T}[f(\xi_{x\delta})-f(\xi)]+\sum_{x:x\in T}\sum_{y:y\thicksim x}\lambda\rho(x)\rho(y)[f(\xi_{xy})-f(\xi)]
\]
for $f\in D(\mathcal{A})$, where
\[\xi_{x\delta}(z)=\left\{ \begin{array}{ll}
\xi(x),& \ {\rm if}\ z\neq x,\\
0,& \ {\rm if}\ z=x.
\end{array} \right.~{\rm and}~~~
\xi_{xy}(z)=\left\{ \begin{array}{ll}
\xi(x),& \ {\rm if}\ z\neq x,\\
\xi(x)+\xi(y),& \ {\rm if}\ z=x.
\end{array} \right.\]
$D(\mathcal{A})$ is the class of all continuous functions $f$ on $\{0,1,2,\cdots\}^{T}$ which has continuous first partial derivatives $f_x$ with respect to $\xi(x)$ and satisfies some $Lipschitz$ condition. (See section 9.1 of \cite{ips} for more details)

To see the relationship between $\{\xi_{t}\}_{t\geq0}$ and $\{\eta_t\}_{t\geq0}$, consider $\widetilde{\eta}_{t}$ given by $\widetilde{\eta}_{t}(x)=I_{\{{\xi}_{t}(x)\geq1\}}$. Then, $\widetilde{\eta}(x)$ flips from 1 to 0 at rate 1 and  flips from 0 to 1 at rate
$\sum_{y:y\thicksim x}\lambda\rho(x)\rho(y)\widetilde{\eta}(y).$ As a result, $\widetilde{\eta}_{t}={\eta}_{t}$ (in distribution).

We now give two lemmas which are critical to get a lower bound of $\lambda_e.$ The former one helps us convert estimating probability into solving ODE. The latter one is a basic fact about the simple random walk. Recall that, $C_t=\{x\in T:\eta_t(x)=1\}$ is the set of the infected individuals at time $t$ and we assume $C_0=\{O\}.$ Recall also that, $\mathbb{P}_{\lambda}(\cdot)=EP_{\lambda}^{\omega}(\cdot)$ is the annealed measure.
\begin{lemma}
For any $t\geq0,$ $\mathbb{P}_{\lambda}^d(C_t\neq\varnothing)\leq \mathbb{E}_{\lambda}^d\xi_t(O).$
\end{lemma}
\textbf{Proof.} By the self duality of the contact process (we give a proof in \textbf{Appendix} \textit{5.2}), $\forall \omega\in\Omega~{\rm and}~t\geq0,$
\[P_{\lambda}^{\omega}(C_t\neq\varnothing)=P_{\lambda}^{\omega}(\eta_t(O)=1).\]

Consequently,  $\mathbb{P}_{\lambda}^d(C_t\neq\varnothing)=\mathbb{P}_{\lambda}^d(\eta_t(O)=1)$.

Hence, \textbf{Lemma 3.1} follows from the relationship between $\{\xi_{t}\}_{t\geq0}$ and $\{\eta_t\}_{t\geq0}$ and the $Markov's\ Inequality:$
\[~~~~~~~\mathbb{P}_{\lambda}^d(\eta_t(O)=1)=\mathbb{P}_{\lambda}^d(\widetilde{\eta}_t(O)=1)=\mathbb{P}_{\lambda}^d(\xi_t(O)\geq1)\leq \mathbb{E}_{\lambda}^d\xi_t(O).~~~~~~\Box\]

Let $\{X_n\}_{n\geq0}$ be the random walk on $T$ with transition probability: \[p(x,y)=\frac{1}{\deg(x)}\] for $y\thicksim x.$ Suppose this random walk is with initial position $X_0=O.$ The law of $\{X_n\}_{n\geq0}$ is denoted by $\widehat{{P}}$ . The expectation corresponding to $\widehat{{P}}$ is denoted by $\widehat{E}.$ Let $|x|=d_{T}(O,x)$, where $d_{T}(\cdot,\cdot)$ is the graph distance.
\begin{lemma}
$\forall x\in(0,1]$ and $n\geq0$, $\widehat{E}x^{|X_n|}\leq[\frac{dx}{d+1}+\frac{1}{(d+1)x}]^n.$
\end{lemma}
\textbf{Proof.}
Let $\{Z_n\}_{n\geq0}$ be a simple random walk on $\mathbb{Z}$ with initial position $0$ and with transition probability: $p(i,i+1)=1-p(i,i-1)=\frac{d}{d+1}.$
Note $Ex^{Z_n}=[Ex^{Z_1-Z_0}]^n=[\frac{dx}{d+1}+\frac{1}{(d+1)x}]^n.$ Hence, we get \textbf{Lemma 3.2} since we can construct a coupling of $\{X_n\}_{n\geq0}$ and $\{Z_n\}_{n\geq0}$ so that $|X_n|\geq Z_n$ for $n\geq0$.~~~~~~$\Box$
\newline
\newline
  We now give the proof of  \eqref{2.2}.\\
\textbf{Proof of \eqref{2.2}.} For $x\in T$ and given vertex weights $\{\rho(x,\omega),{x\in T}\},$ according to the generator $\mathcal{A}$ and Hille-Yosida Theorem, we have the following linear ODE:
\begin{align}\label{ode0}
\frac{d}{dt}E_{\lambda}^{\omega}\xi_t(x)=-E_{\lambda}^{\omega}\xi_t(x)+\sum_{y:y\thicksim x}\lambda\rho(x,\omega)\rho(y,\omega)E_{\lambda}^{\omega}\xi_t(y).
\end{align}
Here, we give an intuitive interpretation for \eqref{ode0}. The Hille-Yosida Theorem shows the following relationship between the semigroup and generator:
\[\frac{d}{dt}E[f(\xi_t)]=E[\mathcal{A}f(\xi_t)].\]
Now, \eqref{ode0} follows by taking $f(\xi)=\xi(x).$ For more rigorous proof, \text{Theorem 1.27.} in Chapter 9 of \cite{ips} is a good reference. The point of  \text{Theorem 1.27.} is that we can apply the Hille-Yosida Theorem to the general linear systems. We can follow the same strategy of \text{Theorem 1.27.} to prove \eqref{ode0}, so we omit its proof here.

Let $G_{\omega}$ be the $T\times T$ matrix with entries
\[G_{\omega}(x,y)=\left\{ \begin{array}{ll}
\lambda\rho(x,\omega)\rho(y,\omega),& \ x\thicksim y,\\
0,&\rm otherwise.
\end{array} \right.\]
and $I$ denote the $T\times T$ identity matrix. Then, equation \eqref{ode0} can be written as the following linear ODE:
\begin{align}\label{ode}
\frac{d}{dt}E_{\lambda}^{\omega}\xi_t=(G_{\omega}-I)E_{\lambda}^{\omega}\xi_t.
\end{align}

Hence, by the standard theory of ODE,  \eqref{ode} has the unique solution
\[E_{\lambda}^{\omega}\xi_t=e^{-t}e^{tG_{\omega}}\xi_0,\]
where $e^{tG_{\omega}}=\sum_{n\geq0}\frac{t^nG_{\omega}^n}{n!}.$

Note $\xi_0\equiv1$, then,
\begin{align}\label{solution}
E_{\lambda}^{\omega}\xi_t(O)=e^{-t}\sum_{n=0}^{\infty}\sum_{x:x\in T}\frac{t^nG_{\omega}^n(O,x)}{n!}.
\end{align}

For $n\geq1,$ if there exists $\{x_i;0\leq i\leq n\}$ such that $x_0=O$, $x_i\in T,$ and $x_{j+1}\sim x_j$ for $0\leq j\leq n-1,$ then we take  $\overrightarrow{x}=(x_0,x_1,\cdot\cdot\cdot,x_n)$ and say $\overrightarrow{x}$ is a path of  length $n$ starting at $O$ . The set of all of such paths is denoted by $V_n$. Note, a path can backtrack.

Then, by \eqref{solution},
\[E_{\lambda}^{\omega}\xi_t(O)
=e^{-t}\sum_{n=0}^{\infty}\frac{t^n\lambda^n}{n!}\Big(\sum_{\overrightarrow{x}\in V_n}\prod_{j=0}^{n-1}\rho(x_j,\omega)\rho(x_{j+1},\omega)\Big),\]
where $\overrightarrow{x}=(x_0,x_1,\cdots,x_n).$

Therefore,
\begin{align}\label{fso}
\mathbb{E}_{\lambda}^d\xi_t(O)
=e^{-t}\sum_{n=0}^{\infty}\frac{t^n\lambda^n}{n!}\Big(\sum_{\overrightarrow{x}\in V_n}E\prod_{j=0}^{n-1}\rho(x_j,\omega)\rho(x_{j+1},\omega)\Big).
\end{align}

Note, $\forall\overrightarrow{x}=(x_0,x_1,\cdots,x_n)\in V_n,$ let $\overrightarrow{y}=(y_0,y_1,\cdots,y_{|x_n|})\in V_{|x_n|}$, where $y_0=O$ and $ y_{|x_n|}=x_n$. Briefly, the path $\overrightarrow{y}$ is from $O$ to $x_n$ and does not backtrack. Then, the vertices in the path $\overrightarrow{y}$ are contained in the path $\overrightarrow{x}$ by the structure of the tree. Set $\widetilde{\rho}=\frac{\rho}{M}\leq1$ and note $\overrightarrow{y}$ consists of $|x_n|+1$ different vertices,
\begin{align}\label{fso1}
E\prod_{j=0}^{n-1}\rho(x_j,\omega)\rho(x_{j+1},\omega)
&=M^{2n}E\prod_{j=0}^{n-1}\widetilde{\rho}(x_j,\omega)\widetilde{\rho}(x_{j+1},\omega)\notag\\
&\leq M^{2n}E\prod_{j=0}^{|x_n|-1}\widetilde{\rho}(y_j,\omega)\widetilde{\rho}(y_{j+1},\omega)\notag\\
&\leq M^{2n}[E{\widetilde{\rho}}^2]^{|x_n|-1}.
\end{align}

Hence, by \eqref{fso} and \eqref{fso1},
\begin{align}\label{fso2}
\mathbb{E}_{\lambda}^d\xi_t(O)
&\leq e^{-t}\sum_{n=0}^{\infty}\frac{t^n\lambda^n}{n!}\big[\sum_{\overrightarrow{x}\in V_n}(E{\widetilde{\rho}}^2)^{|x_n|-1}\big]M^{2n}\notag\\
&=e^{-t}\sum_{n=0}^{\infty}\frac{t^n\lambda^nM^{2n}}{n!E{\widetilde{\rho}}^2}\big[\sum_{\overrightarrow{x}\in V_n}(E{\widetilde{\rho}}^2)^{|x_n|}\big]
\end{align}

Note, $\forall x\in T^d,$ the degree of $x$ satisfies $\deg(x)\leq d+1$, then
\begin{align}\label{k1}
\sum_{\overrightarrow{x}\in V_n}(E{\widetilde{\rho}}^2)^{|x_n|}
\leq (d+1)^n\sum_{\overrightarrow{x}\in V_n}\Big(\prod_{j=0}^{n-1}\frac{1}{\deg(x_j)}\Big)\Big(E{\widetilde{\rho}}^2\Big)^{|x_n|}
\end{align}
Recall the definition of the random walk $\{X_n\}_{n\geq1},$ we have,
\[\widehat{{P}}(X_j=x_j,0\leq j\leq n)=\prod_{j=0}^{n-1}\frac{1}{\deg(x_j)},~\forall\overrightarrow{x}=(x_0,x_1,\cdots,x_n)\in V_n.\]
Therefore,
\begin{align}\label{k2}
\sum_{\overrightarrow{x}\in V_n}\Big(\prod_{j=0}^{n-1}\frac{1}{\deg(x_j)}\Big)\Big(E{\widetilde{\rho}}^2\Big)^{|x_n|}
=\widehat{{E}}[(E{\widetilde{\rho}}^2)^{|X_n|}]
\end{align}
So, by \eqref{k1} and \eqref{k2},
\begin{align}\label{k3}
\sum_{\overrightarrow{x}\in V_n}(E{\widetilde{\rho}}^2)^{|x_n|}\leq (d+1)^n\widehat{{E}}[(E{\widetilde{\rho}}^2)^{|X_n|}].
\end{align}
Hence, by \eqref{fso2} and \eqref{k3},
\begin{align}\label{be}
\mathbb{E}_{\lambda}^d\xi_t(O)\leq e^{-t}\sum_{n=0}^{\infty}\frac{t^n\lambda^nM^{2n}(d+1)^n}{n!E{\widetilde{\rho}}^2}{\widehat{E}}[(E{\widetilde{\rho}}^2)^{|X_n|}]
\end{align}
Since $E({\widetilde{\rho}}^2)\in(0,1],$ then, by \textbf{Lemma 3.2},
\begin{align}\label{be1}
\widehat{{E}}[(E{\widetilde{\rho}}^2)^{|X_n|}]\leq [\frac{dE({\widetilde{\rho}}^2)}{d+1}+\frac{1}{(d+1)E({\widetilde{\rho}}^2)}]^n.
\end{align}
Hence, by \eqref{be} and \eqref{be1},
\begin{align}\label{ee}
\mathbb{E}_{\lambda}^d\xi_t(O)
&\leq e^{-t}\sum_{n=0}^{\infty}\frac{t^n\lambda^nM^{2n}(d+1)^n}{n!E{\widetilde{\rho}}^2}[\frac{dE({\widetilde{\rho}}^2)}{d+1}+\frac{1}{(d+1)E({\widetilde{\rho}}^2)}]^n\notag\\
&=\big(E({\widetilde{\rho}}^2)\big)^{-1}\exp\big\{t\big[\lambda M^2\big(dE({\widetilde{\rho}}^2)+\frac{1}{E({\widetilde{\rho}}^2)}\big)-1\big]\big\}.
\end{align}
Therefore, by \textbf{Lemma 3.1} and \eqref{ee},
\[\mathbb{P}_{\lambda}^d(C_t^O\neq\varnothing)\leq
\big({E}({\widetilde{\rho}}^2)\big)^{-1}\exp\big\{t\big[\lambda M^2\big(dE({\widetilde{\rho}}^2)+\frac{1}{E({\widetilde{\rho}}^2)}\big)-1\big]\big\}.\]
As a result, for $\lambda<(dE{\rho}^2+\frac{M^4}{E{\rho}^2})^{-1},$
\begin{align}
\limsup_{t\rightarrow\infty}\frac{1}{t}\log \mathbb{P}_{\lambda}^d(C_t^O\neq\varnothing)
&\leq\lambda M^2\big(dE({\widetilde{\rho}}^2)+\frac{1}{E({\widetilde{\rho}}^2)}\big)-1\notag\\
&=\lambda(dE{\rho}^2+\frac{M^4}{E{\rho}^2})-1<0\notag
\end{align}
Thus, $\lambda_c\geq\lambda_e\geq(dE{\rho}^2+\frac{M^4}{E{\rho}^2})^{-1}>0.$~~~~~~$\Box$

\section{Upper bound for $\lambda_c$}
\label{S:4}
In this section, we use the second moments method to get an upper bound of $\lambda_c$.
We first consider the SIR model, which can dominate the process $\eta_t$ from below. Recall that, $\mathbb{P}_{\lambda}(\cdot)=EP_{\lambda}^{\omega}(\cdot)$ is the annealed measure and $C_t=\{x\in T:\eta_t(x)=1\}$ is the set of the infected individuals at time $t$ and we assume $C_0=\{O\}$.

Let $\{\zeta_t\}_{t\geq0}$ be Markov process on $\{-1,0,1\}^T$. $\forall t\geq0,$ let $S_t=\{x\in T:\zeta_t(x)=0\}$, $I_t=\{x\in T:\zeta_t(x)=1\},$ $R_t=\{x\in T:\zeta_t(x)=-1\}.$ Assume that, $\zeta_0(O)=1$ and $\zeta_0(x)=0$ for any $x\neq O$ when $t=0.$ When the vertex weights $\{\rho(x,\omega),x\in T\}$ is given, $\{\zeta_t\}_{t\geq0}$ flips only in the following two situations:
If $x\in I_t$, then $\zeta_t(x)$ flips to $-1$ at rate $1.$
If $y$ is a son of $x$, and $x\in I_t,$ $y\in S_t,$ then $\zeta_t(y)$ flips to $1$ at rate $\lambda\rho(x,\omega)\rho(y,\omega).$ Individuals in $R_t$ will always be in $R_t$. In brief, $\zeta_t$ is a non-backtrack contact process on $T$ from the root $O$.

This model is a generalize of the classic contact model. $R_t$ can be interpreted as the removed or quarantined individuals, so they slow down the progression of the epidemic. Hence, $I_t\subseteq C_t$  for any $t\geq0$ in the sense of coupling (see section 2.1 of \cite{ips} for more details) when $\{\eta_t\}_{t\geq0}$ and $\{\zeta_t\}_{t\geq0}$ have the same infection rate $\lambda$ and
vertex weights $\{\rho(x,\omega),{x\in T}\}.$ Hence, $\forall\omega\in\Omega$,
\[P_{\lambda}^{\omega}(C_t\neq\varnothing, \forall t\geq0)\geq P_{\lambda}^{\omega}(I_t\neq\varnothing, \forall t\geq0),\]
 and then, $\forall\lambda>0$ and $d\geq1,$
\[\mathbb{P}_{\lambda}^{d}(C_t\neq\varnothing, \forall t\geq0)\geq \mathbb{P}_{\lambda}^{d}(I_t\neq\varnothing, \forall t\geq0).\]

Let $I_{\infty}=\bigcup_{t\geq0}I_t$. Then ,
\[\{I_t\neq\varnothing, \forall t\geq0\}=\{|I_{\infty}|=+\infty\}.\]

Thus,
\[\mathbb{P}_{\lambda}^{d}(C_t\neq\varnothing, \forall t\geq0)\geq \mathbb{P}_{\lambda}^{d}(|I_{\infty}|=+\infty).\]

Let $L_n=\{x\in T:|x|=n,x\in I_{\infty}\}$, then
\[\mathbb{P}_{\lambda}^{d}(|I_{\infty}|=+\infty)\geq \mathbb{P}_{\lambda}^{d}(\bigcap_{n=0}^{\infty}\{|L_n|>0\})=\lim_{n\rightarrow\infty}\mathbb{P}_{\lambda}^{d}(|L_n|>0).\]

Since $|L_n|$ is a nonnegative random variable, so by $Cauchy's\ Inequality$:
\[\mathbb{E}_{\lambda}^d|L_n|=\mathbb{E}_{\lambda}^d|L_n|I_{\{|L_n|>0\}}\leq (\mathbb{E}_{\lambda}^d|L_n|^2)^{\frac{1}{2}}(\mathbb{P}_{\lambda}^d(|L_n|>0))^{\frac{1}{2}},\]

Hence, by all of the above,
\begin{align}\label{17}
\mathbb{P}_{\lambda}^{d}(C_t\neq\varnothing, \forall t\geq0)\geq \lim_{n\rightarrow\infty}\mathbb{P}_{\lambda}^d(|L_n|>0)\geq\limsup_{n\rightarrow\infty}\frac{(\mathbb{E}_{\lambda}^d|L_n|)^2}{\mathbb{E}_{\lambda}^d|L_n|^2}.
\end{align}

To handle with $\mathbb{E}_{\lambda}^d|L_n|$ and $\mathbb{E}_{\lambda}^d|L_n|^2$, we first introduce the non-backtrack random walk $\{S_n\}_{n\geq0}$ on $T$. Let $S_0=O$, and when $S_n=x$ and $y$ is any son of $x$,
\[P(S_{n+1}=y|S_n=x)=\frac{1}{d}.\]

Let $\{\widehat{S}_n\}_{n\geq0}$ be an independent copy of $\{S_n\}_{n\geq0}$ and both of them are defined on $(\widetilde{{\Omega}},\widetilde{{\mathcal{F}}},\widetilde{{P}}).$

Now consider the first moment $\mathbb{E}_{\lambda}^d|L_n|$ and write it as the following form:
\begin{align}\label{18}
\mathbb{E}_{\lambda}^d|L_n|=\mathbb{E}_{\lambda}^d\sum_{x:|x|=n}I_{\{x\in I_{\infty}\}}=\sum_{x:|x|=n}\mathbb{P}_{\lambda}^d(x\in I_{\infty}).
\end{align}

Now, fix some $x\in T$, $|x|=n$ and consider $\mathbb{P}_{\lambda}^d(x\in I_{\infty}).$ Let $\{H_x\}_{x\in T}$ be $i.i.d$ exponential times with rate $1.$ If $y$ is a son of $x$ (take $x\rightarrow y$ as briefly), let $U_{x,y}$ be
exponential time with rate $\lambda\rho(x,\omega)\rho(y,\omega)$ (or $\lambda\rho(x)\rho(y)$ for simplicity). Assume all of the exponential times are independent.
By the structure of the tree, there is a unique path $\{O=x_0\rightarrow x_1\rightarrow\cdots\rightarrow x_n=x\}$ such that $|x_n|=n.$
Therefore,
\begin{align}
P_{\lambda}^{\omega}(x\in I_{\infty})
&=P_{\lambda}^{\omega}(U_{x_i,x_{i+1}}<H_{x_i},\forall i\in[0,n-1])\notag\\
&=\prod_{i=0}^{n-1}[\frac{\lambda\rho(x_i,\omega)\rho(x_{i+1},\omega)}{1+\lambda\rho(x_i,\omega)\rho(x_{i+1},\omega)}].\notag
\end{align}

Hence,
\begin{align}\label{19}
\mathbb{P}_{\lambda}^{d}(x\in I_{\infty})
=E\prod_{i=0}^{n-1}[\frac{\lambda\rho(x_i,\omega)\rho(x_{i+1},\omega)}{1+\lambda\rho(x_i,\omega)\rho(x_{i+1},\omega)}]>0.
\end{align}

Since $\{\rho(x),{x\in T}\}$ are $i.i.d.$ random variables, then $\forall y\in T$ and $|y|=n,$
\[\mathbb{P}_{\lambda}^{d}(x\in I_{\infty})=\mathbb{P}_{\lambda}^{d}(y\in I_{\infty}).\]

Thus, by equation \eqref{18} and \eqref{19},
\begin{align}\label{20}
\mathbb{E}_{\lambda}^d|L_n|=d^nE\prod_{i=0}^{n-1}[\frac{\lambda\rho(x_i,\omega)\rho(x_{i+1},\omega)}{1+\lambda\rho(x_i,\omega)\rho(x_{i+1},\omega)}].
\end{align}

Hence, by the definition of the non-backtrack random walk $\{S_n\}_{n\geq0}$ and \eqref{20},
\begin{align}\label{1mom}
\mathbb{E}_{\lambda}^d|L_n|=d^nE\prod_{i=0}^{n-1}[\frac{\lambda\rho(S_i,\omega)\rho(S_{i+1},\omega)}{1+\lambda\rho(S_i,\omega)\rho(S_{i+1},\omega)}].
\end{align}

Note, in the right side of equation \eqref{1mom}, we do not take the expectation with respect to $\{S_n\}_{n\geq0}$. Hence, it is a random variable on $(\widetilde{{\Omega}},\widetilde{{\mathcal{F}}},\widetilde{{P}})$ and a constant almost everywhere with respect to $\widetilde{{P}}.$

We now turn to calculate $\mathbb{E}_{\lambda}^d|L_n|^2.$ Since,
 \[|L_n|=\sum_{x:|x|=n}I_{\{x\in L_n\}},\]

Then,
\begin{align}\label{22}
\mathbb{E}_{\lambda}^d|L_n|^2=\sum_{x:|x|=n}\sum_{y:|y|=n}\mathbb{P}_{\lambda}^d(x\in L_n,y\in L_n).
\end{align}

For convenience, for any $x\rightarrow u$ and $y\rightarrow v$, let $F(x,y;u,v)=P_{\lambda}^{\omega}(U_{x,u}<H_x,U_{y,v}<H_{y}).$ Then, we have the following result by direct calculation.
\[F(x,y;u,v)
\left\{ \begin{array}{llll}
=\frac{{\lambda}^2\rho(x)\rho(y)\rho(u)\rho(v)}{[1+\lambda\rho(x)\rho(u)][1+\lambda\rho(y)\rho(v)]},& \ {\rm if}\ x\neq y,\\
=\frac{\lambda\rho(x)\rho(u)}{1+\lambda\rho(x)\rho(u)},& \ {\rm{if}}\ x=y\ {\rm and}\ u=v ,\\
\leq\frac{2\lambda^2\rho^2(x)\rho(u)\rho(v)}{[1+\lambda\rho(x)\rho(u)]
[1+\lambda\rho(x)\rho(v)]},& \ {\rm if}\ x=y\ {\rm and}\ u\neq v,\\
=0,& \ {\rm otherwise}.
\end{array} \right.\]

Now, fix $x\in T$ and $|x|=n$, $y\in T$ and $|y|=n$. Take $x\wedge y$ as the common ancestor furthest from $O$, and \\
$O=x_0\rightarrow x_1\rightarrow\cdots\rightarrow x_{k-1}\rightarrow x_k=x\wedge y\rightarrow z=x_{k+1}\rightarrow\cdots\rightarrow x_n=x,$\\
$O=x_0=y_0\rightarrow x_1=y_1\rightarrow\cdots\rightarrow x_{k-1}=y_{k-1}\rightarrow x_k=y_k=x\wedge y\rightarrow w=y_{k+1}\rightarrow\cdots\rightarrow y_n=y.$
By the independence of the exponential times,
\begin{align}\label{23}
\mathbb{P}_{\lambda}^d(x\in L_n,y\in L_n)
&=EP_{\lambda}^{\omega}(\forall0\leq i\leq n-1, U_{x_i,x_{i+1}}<H_{x_i},U_{y_i,y_{i+1}}<H_{y_i})\notag\\
&=E[\prod_{i=0}^{n-1}F(x_i,y_i;x_{i+1},y_{i+1})].
\end{align}

Take $\Lambda_n$ as the set of all non-backtrack paths with length $n,$ then, by \eqref{22} and \eqref{23},
\begin{align}\label{equ 1}
\mathbb{E}_{\lambda}^d|L_n|^2
&=\sum_{\overrightarrow{x}\in\Lambda_n}\sum_{\overrightarrow{y}\in\Lambda_n}E[\prod_{i=0}^{n-1}F(x_i,y_i;x_{i+1},y_{i+1})]
\end{align}
where $\overrightarrow{x}=(x_0,x_1,\cdots,x_n)\in\Lambda_n$ and $\overrightarrow{y}=(y_0,y_1,\cdots,y_n)\in\Lambda_n$.

On the other hand, we can take $\overrightarrow{x}$ and $\overrightarrow{y}$
as an sample path of $S_n$ and $\widehat{S}_n$ respectively. Note $\widetilde{P}(S_n=\overrightarrow{x})=\widetilde{P}(\widehat{S}_n=\overrightarrow{y})=\frac{1}{d^n}.$

Therefore,
\begin{align}\label{25}
&(\widetilde{E}\times E)[\prod_{i=0}^{n-1}F(S_i,\widehat{{S}}_i;S_{i+1},\widehat{S}_{i+1})]\notag\\
=&\sum_{\overrightarrow{x}\in\Lambda_n}\sum_{\overrightarrow{y}\in\Lambda_n}\frac{1}{d^{2n}}E[\prod_{i=0}^{n-1}F(x_i,y_i;x_{i+1},y_{i+1})].
\end{align}

Thus, by \eqref{equ 1} and (25),
\begin{align}\label{26}
\mathbb{E}_{\lambda}^d|L_n|^2=d^{2n}(\widetilde{E}\times E)[\prod_{i=0}^{n-1}F(S_i,\widehat{{S}}_i;S_{i+1},\widehat{S}_{i+1})]
\end{align}

As a consequence of \eqref{17}, \eqref{1mom} and \eqref{26}, in order that the process survives, it suffices to show that
\begin{align}\label{condi}
\limsup_{n\rightarrow\infty}\widetilde{E}\big[\frac{E\prod_{i=0}^{n-1}F(S_i,\widehat{{S}}_i;S_{i+1},\widehat{S}_{i+1})}
{(E\prod_{i=0}^{n-1}\frac{\lambda\rho(S_i)\rho(S_{i+1})}{1+\lambda\rho(S_i)\rho(S_{i+1})})^2}\big]<+\infty.
\end{align}

In order to simplify  \eqref{condi}, take
\[f_n=\frac{E\prod_{i=0}^{n-1}F(S_i,\widehat{{S}}_i;S_{i+1},\widehat{S}_{i+1})}
{(E\prod_{i=0}^{n-1}\frac{\lambda\rho(S_i)\rho(S_{i+1})}{1+\lambda\rho(S_i)\rho(S_{i+1})})^2}.\]

Note $f_n$ is a function of $\{S_i,\widehat{S}_{i}:0\leq i\leq n\}$. So, it is a random variable on $(\widetilde{{\Omega}},\widetilde{{\mathcal{F}}},\widetilde{{P}})$. We will use the notation $f_n(\tilde{\omega})$ when emphasizing the randomness of $f_n.$   Hence, \eqref{condi} is equivalent to
\[\limsup_{n\rightarrow\infty}\widetilde{E}f_n(\tilde{\omega})<+\infty.\]

We summarize the above results in the following Proposition.
\begin{proposition}
If~$\limsup_{n\rightarrow\infty}\widetilde{E}f_n(\tilde{\omega})<+\infty$ for some $\lambda\in(0,\infty),$ then
$\mathbb{P}_{\lambda}^{d}(C_t^O\neq\varnothing, \forall t\geq0)>0.$ Hence, $\lambda_c\leq\lambda.$~~~~~~$\Box$
\end{proposition}
Now we give a proposition with which we can get the proof of \textbf{Theorem 2.1}.
\begin{proposition}
If $\frac{(1+\lambda M^2)^2}{\lambda E\rho^2}<d$ for some $\lambda\in(0,\infty),$ then
\[\limsup_{n\rightarrow\infty}\widetilde{E}f_n(\tilde{\omega})<+\infty.\]
\end{proposition}
\textbf{Proof.} Let $\tau=\sup\{n\geq0:S_n=\widehat{S}_n\},$ then
\begin{align}\label{28}
\widetilde{P}(\tau=k)=\frac{1}{d^k}(1-\frac{1}{d}),\ \forall k\geq0.
\end{align}

Thus, by \eqref{28}, for $s\in(0,d),$
\begin{align}\label{29}
\widetilde{E}s^{\tau}=\sum_{k=0}^{\infty}\frac{s^k}{d^k}(1-\frac{1}{d})<\infty.
\end{align}

Write $\widetilde{E}f_n(\tilde{\omega})$ in the following form:
\[\widetilde{E}f_n(\tilde{\omega})
=\widetilde{E}[f_n(\tilde{\omega})I_{\{\tau<n\}}]+\widetilde{E}[f_n(\tilde{\omega})I_{\{\tau\geq n\}}]:=\textrm{I}_n+\textrm{II}_n.\]

We consider $\textrm{II}_n$ first and note that on the event $\{\tau\geq n\},$
\begin{align}\label{30}
E\prod_{i=0}^{n-1}F(S_i,\widehat{{S}}_i;S_{i+1},\widehat{S}_{i+1})
=E\prod_{i=0}^{n-1}\frac{\lambda\rho(S_i)\rho(S_{i+1})}{1+\lambda\rho(S_i)\rho(S_{i+1})}.
\end{align}

Hence, by \eqref{30}, on the event $\{\tau\geq n\}$,
\begin{align}\label{31}
f_n(\tilde{\omega})&=\big{(}E\prod_{i=0}^{n-1}\frac{\lambda\rho(S_i)\rho(S_{i+1})}{1+\lambda\rho(S_i)\rho(S_{i+1})}\big{)}^{-1}\notag\\
&\leq\big(E\prod_{i=0}^{n-1}\frac{\lambda\rho(S_i)\rho(S_{i+1})}{1+\lambda M^2}\big)^{-1}\notag\\
&=[\frac{(1+\lambda M^2)}{\lambda E\rho^2}]^n\frac{E\rho^2}{(E\rho)^2}.
\end{align}

Thus, according to \eqref{31},
\begin{align}\label{32}
\textrm{II}_n\leq[\frac{(1+\lambda M^2)}{\lambda E\rho^2}]^n\frac{E\rho^2}{(E\rho)^2}P(\tau\geq n)=
[\frac{(1+\lambda M^2)}{d\lambda E\rho^2}]^n\frac{E\rho^2}{(E\rho)^2}.
\end{align}

Hence, by \eqref{32}, $\textrm{II}_n$ is bounded when
\begin{align}\label{33}
\frac{(1+\lambda M^2)}{d\lambda E\rho^2}\leq1.
\end{align}

Now, we turn to handle with $\textrm{I}_n$ and notice on  $\{\tau<n\}$,
\begin{align}\label{34}
f_n(\tilde{\omega})=\frac{E[\prod_{i=0}^{\tau}F(S_i,\widehat{{S}}_i;S_{i+1},\widehat{S}_{i+1})]
E[\prod_{i=\tau+1}^{n-1}F(S_i,\widehat{{S}}_i;S_{i+1},\widehat{S}_{i+1})]}
{\{E\prod_{i=0}^{\tau}[\frac{\lambda\rho(S_i)\rho(S_{i+1})}{1+\lambda\rho(S_i)\rho(S_{i+1})}]\}^2
\{E\prod_{i=\tau+1}^{n-1}[\frac{\lambda\rho(S_i)\rho(S_{i+1})}{1+\lambda\rho(S_i)\rho(S_{i+1})}]\}^2}.
\end{align}

(The product over the empty set is taken to be $1.$)

Since $\{\rho(x),x\in T\}$ is an $i.i.d.$ sequence, then,
\begin{align}\label{35}
E\prod_{i=\tau+1}^{n-1}[\frac{\lambda\rho(S_i)\rho(S_{i+1})}{1+\lambda\rho(S_i)\rho(S_{i+1})}]
=E\prod_{i=\tau+1}^{n-1}[\frac{\lambda\rho(\widehat{S}_i)\rho(\widehat{S}_{i+1})}{1+\lambda\rho(\widehat{S}_i)\rho(\widehat{S}_{i+1})}],
\end{align}
\begin{align}\label{36}
&E[\prod_{i=\tau+1}^{n-1}F(S_i,\widehat{{S}}_i;S_{i+1},\widehat{S}_{i+1})]\notag\\
=&E\prod_{i=\tau+1}^{n-1}[\frac{\lambda\rho(S_i)\rho(S_{i+1})}{1+\lambda\rho(S_i)\rho(S_{i+1})}]
E\prod_{i=\tau+1}^{n-1}[\frac{\lambda\rho(\widehat{S}_i)\rho(\widehat{S}_{i+1})}{1+\lambda\rho(\widehat{S}_i)\rho(\widehat{S}_{i+1})}].
\end{align}

So by \eqref{34}, \eqref{35} and \eqref{36}, on $\{\tau<n\}$,
\[f_n(\tilde{\omega})=\frac{E[\prod_{i=0}^{\tau}F(S_i,\widehat{{S}}_i;S_{i+1},\widehat{S}_{i+1})]}
{\{E\prod_{i=0}^{\tau}[\frac{\lambda\rho(S_i)\rho(S_{i+1})}{1+\lambda\rho(S_i)\rho(S_{i+1})}]\}^2}.\]

Denote $E[\prod_{i=0}^{\tau}F(S_i,\widehat{{S}}_i;S_{i+1},\widehat{S}_{i+1})]$ by $\textrm{I}_{n,1}$ and
$\{E\prod_{i=0}^{\tau}[\frac{\lambda\rho(S_i)\rho(S_{i+1})}{1+\lambda\rho(S_i)\rho(S_{i+1})}]\}^2$ by $\textrm{I}_{n,2}.$ Then,
\begin{align}\label{37}
\textrm{I}_n=\widetilde{E}\big{[}\frac{\textrm{I}_{n,1}}{\textrm{I}_{n,2}}I_{\{\tau<n\}}\big{]}.
\end{align}

For $\textrm{I}_{n,1}$, note that $S_{\tau}=\widehat{{S}}_{\tau}$ and $S_{\tau+1}\neq\widehat{S}_{\tau+1}$, then
\begin{align}\label{38}
F(S_{\tau},\widehat{{S}}_{\tau};S_{\tau+1},\widehat{S}_{\tau+1})
&\leq\frac{2{\lambda}^2{\rho^2(S_{\tau})}\rho(S_{\tau+1})\rho(\widehat{S}_{\tau+1})}
{[1+\lambda\rho(S_{\tau})\rho(S_{\tau+1})][1+\lambda\rho(\widehat{S}_{\tau})\rho(\widehat{S}_{\tau+1})]}\notag\\
&\leq2{\lambda}^2M^2\rho(S_{\tau+1})\rho(\widehat{S}_{\tau+1}).
\end{align}

Hence, by \eqref{38} and the independence of $\{\rho(x),{x\in T}\}$,
\begin{align}\label{39}
\textrm{I}_{n,1}&=E[F(S_{\tau},\widehat{{S}}_{\tau};S_{\tau+1},\widehat{S}_{\tau+1})\prod_{i=0}^{\tau-1}F(S_i,\widehat{{S}}_i;S_{i+1},\widehat{S}_{i+1})]\notag\\
&\leq E\big[\prod_{i=0}^{\tau-1}\frac{\lambda\rho(S_i)\rho(S_{i+1})}{1+\lambda\rho(S_i)\rho(S_{i+1})}\big(2{\lambda}^2M^2\rho(S_{\tau+1})\rho(\widehat{S}_{\tau+1})\big)\big]\notag\\
&\leq2{\lambda}^2M^2E[\prod_{i=0}^{\tau-1}\lambda\rho(S_i)\rho(S_{i+1})]\{E\rho(S_{\tau+1})\}^2\notag\\
&=2M^2(E\rho)^4{\lambda}^{2+\tau}(E\rho^2)^{\tau-1}.
\end{align}

Therefore, by \eqref{39}, we get the upper bound of $\textrm{I}_{n,1}$:
\begin{align}\label{40}
\textrm{I}_{n,1}\leq2M^2(E\rho)^4{\lambda}^{2+\tau}(E\rho^2)^{\tau-1}.
\end{align}

For $\textrm{I}_{n,2}$,
\begin{align}\label{41}
\textrm{I}_{n,2}&=\{E\prod_{i=0}^{\tau}[\frac{\lambda\rho(S_i)\rho(S_{i+1})}{1+\lambda\rho(S_i)\rho(S_{i+1})}]\}^2\notag\\
&\geq \{E\prod_{i=0}^{\tau}[\frac{\lambda\rho(S_i)\rho(S_{i+1})}{1+\lambda M^2}]\}^2\notag\\
&=(E\rho)^4(\frac{\lambda}{1+\lambda M^2})^{2\tau+2}(E\rho^2)^{2\tau}.
\end{align}

Thus, by \eqref{40} and \eqref{41}, on $\{\tau<n\}$,
\begin{align}\label{42}
\frac{\textrm{I}_{n,1}}{\textrm{I}_{n,2}}\leq C(\lambda,M)[\frac{(1+\lambda M^2)^2}{\lambda E\rho^2}]^{\tau},
\end{align}
where $C(\lambda,M)=2M^2(1+\lambda M^2)^2(E\rho^2)^{-1}>0.$

Hence, by \eqref{37} and \eqref{42},
\begin{align}\label{43}
\textrm{I}_n&\leq C(\lambda,M)\widetilde{E}\big{\{}[\frac{(1+\lambda M^2)^2}{\lambda E\rho^2}]^{\tau}I_{\{\tau<n\}}\big{\}}\notag\\
&\leq C(\lambda,M)\widetilde{E}\big{\{}[\frac{(1+\lambda M^2)^2}{\lambda E\rho^2}]^{\tau}\big{\}}.
\end{align}

Thus,  by \eqref{29} and \eqref{43}, $\textrm{I}_n$ is bounded when
\begin{align}\label{upcondition}
\frac{(1+\lambda M^2)^2}{d\lambda E\rho^2}<1.
\end{align}

Hence, \textbf{Proposition 4.2} follows from \eqref{33} and \eqref{upcondition}.~~~~~~$\Box$

In section 3, we have shown \eqref{2.2}. \textbf{Proposition 4.1} and \textbf{Proposition 4.2} show that, if $\frac{(1+\lambda M^2)^2}{\lambda E\rho^2}<d$ for some $\lambda\in(0,\infty)$, then, $\lambda_e\leq\lambda_c\leq\lambda<\infty.$ Hence, we accomplish  the proof of \textbf{Theorem 2.1}.

\textbf{Acknowledgments.} The authors are grateful to the financial support from the National Natural Science Foundation of China with grant number 11371040, 11531001 and 11501542.

\section{Appendix}
\subsection{The process is a spin system}
To see $\{\eta_t\}_{t\geq0}$ is a spin system, note that
\[c(x,\eta)\leq \max\{1,\lambda M^2(d+1)\}.\]

Thus $c(x,\eta)$ is a uniformly bounded nonnegative function.
For each $x\in T,$ $\xi(u)=\eta(u)$ for $u\in \{x\}\bigcup\{y:y\sim x\}$ when $\xi$ is close enough to $\eta$ by the topology of $\{0,1\}^{T}.$

Hence $c(x,\xi)=c(x,\eta)$ for all $\xi$ close enough to $\eta$. So $c(x,\eta)$ is continuous in $\eta$ for each $x.$

Then consider the following supremum:
\[\sup_{x\in T}\sum_{u\in T}\sup_{\eta\in X}|c(x,\eta)-c(x,\eta_u)|,\] where $\eta_u\in X$ is defined by
\[\eta_u(v)=\left\{ \begin{array}{ll}
\eta(v),& \ {\rm if}\ v\neq u,\\
1-\eta(v),& \ {\rm if}\ v=u.
\end{array} \right.\]
Hence, $\forall x \in T:$

If $u\sim x,$
\begin{align}\label{45}
\sup_{\eta\in X}|c(x,\eta)-c(x,\eta_u)|=\lambda\rho(x)\rho(u)\leq \lambda M^2,
\end{align}
and the supremum is taken when $\eta(x)=0.$

If $u=x,$
\begin{align}\label{46}
\sup_{\eta\in X}|c(x,\eta)-c(x,\eta_u)|\leq \max\{1,\lambda M^2(d+1)-1\},
\end{align}

If $u\neq x$ and is not a neighbor of $x$, it is obvious that
\begin{align}\label{47}
\sup_{\eta\in X}|c(x,\eta)-c(x,\eta_u)|=\sup_{\eta\in X}0=0.
\end{align}
Hence by \eqref{45}, \eqref{46} and \eqref{47},
\[\sup_{x\in T}\sum_{u\in T}\sup_{\eta\in X}|c(x,\eta)-c(x,\eta_u)|<\infty.\]
Then according to Chapter 3 of \cite{ips}, $\{\eta_t\}_{t\geq0}$ is a spin system.$~~~~~~~~~~~~~~~~~~\Box$
\subsection{Duality for the process}
$\forall\omega\in\Omega$ and $t\geq0$, we have the following self duality for our contact process $\{\eta_t\}_{t\geq0}$ with initial state $\eta_0\equiv1$:
\begin{align}\label{48}
P_{\lambda}^{\omega}(C_t\neq\varnothing)=P_{\lambda}^{\omega}(\eta_t(O)=1),
\end{align}
where $C_t=\{x\in T:\eta_t(x)=1\}$ is the set of the infected individuals at time $t$ and we assume $C_0=\{O\}.$

\textbf{Proof.} In this model, $\{\eta_t\}_{t\geq0}$ has the following transition measure :
\[c(x,\eta)=\left\{ \begin{array}{ll}
1,& \ {\rm if}\ \eta(x)=1,\\
\lambda\sum_{y:y\sim x}\rho(x)\rho(y)\eta(y),& \ {\rm if}\ \eta(x)=0.
\end{array} \right.\]

Let
\begin{align}
&c(x)=1+\lambda\sum_{u:u\sim x}\rho(x)\rho(u),\notag\\
&p(x,\varnothing)=\frac{1}{1+\lambda\sum_{u:u\sim x}\rho(x)\rho(u)},\notag\\
&p(x,\{x,y\})=\frac{\lambda\rho(x)\rho(y)}{1+\lambda\sum_{u:u\sim x}\rho(x)\rho(u)},\ {\rm when}\ y\sim x,\notag\\
&p(x,A)=0,\ \rm otherwise,\notag
\end{align}

and $Y=\{A:A\subset T,|A|<\infty\}$.

Define the dual function \[H(\eta,A)=\Pi_{x\in A}[1-\eta(x)].\]

(The product over the empty set is taken to be 1.)

Then the rates $c(x,\eta)$ can be written in the form:
\[c(x,\eta)=c(x)\{[1-\eta(x)]+[2\eta(x)-1]\sum_{A\in Y}p(x,A)H(\eta,A)\}.\]

Hence, by Chapter 3 of \cite{ips} there is a coalescing dual process $\{A_t\}_{t\geq0}$ on $Y$ with generator:
\[\mathcal{A}f(A)=\sum_{x:x\in A}[f(A\backslash x)-f(A)]+\sum_{x:x\in A}\sum_{y:y\thicksim x}\lambda\rho(x)\rho(y)[f(A\cup\{y\})-f(A)]\]
for $f\in C(Y),$ then we know the dual chain $\{A_t\}_{t\geq0}$ is exactly the contact process $\{C_t\}_{t\geq0}$ itself.

Thus, we get the following duality by \textbf{Theorem 4.13.} in Chapter 3 of \cite{ips}:
For every $\eta\in X,$ $A\in Y$ and $t\geq0$,
\begin{align}\label{dual}
S(t)H(\cdot,A)(\eta)=E^{A}H(\eta,A_t),
\end{align}
where $S(t)$ is the Markov Semigroup of the process $\{\eta_t\}_{t\geq0}$.

Now take $\eta=\eta_0$ and $A=\{O\}$ in equation \eqref{dual}, we get instantly the dual equality \eqref{48} by the definition of the dual function $H.$
\newline\\
{}
\end{document}